\documentclass[twocolumn]{autart}    
\usepackage{amsmath,amssymb}

\usepackage{accents}
\usepackage[toc,page]{appendix}
\usepackage{enumerate}

\newcommand{\bv}[1]{{\mathbf #1}} 
\newcommand{\bt}[1]{{\mathbf{#1}}}

\newcommand{\rom}[1]{{\hat{#1}}}
\newcommand{\mmark}[1]{{\mathring{#1}}}
\newcommand{\R}{\mathbb{R}}

\newcommand{\ct}[1]{{\mathcal{#1}}} 
\newcommand{\cvv}[1]{{\mathfrak{#1}}}

\newcommand{\mexp}[1]{{e^{#1}}}

\newcommand{\Ob}{\bt{Q}}
\newcommand{\Obc}{\ct{Q}}

\newcommand{\acro}[1]{{\texttt{#1}}}  

\usepackage{graphicx}          

\newcommand{\funSpace}[1]{\mathcal{#1}} 
\usepackage{makecell}
\usepackage{multirow}
\usepackage[normalem]{ulem}
\newcommand{\sepsymb}[1]{\dashuline{#1}}
\usepackage{marvosym} 

\begin{document}

\begin{frontmatter}

\title{Effective error estimation for model reduction with inhomogeneous initial conditions\thanksref{footnoteinfo}} 

\thanks[footnoteinfo]{
Corresponding author B.~Liljegren-Sailer. Tel. +49 651 201-3468.
}

\author[Trier]{Bj{\"o}rn Liljegren-Sailer}\ead{bjoern.sailer@uni-trier.de} 

\address[Trier]{Universit{\"a}t Trier, FB IV - Mathematik, Lehrstuhl Modellierung und Numerik, D-54286 Trier, Germany}  

%


\begin{abstract}                          
A priori error bounds have been derived for different balancing-related model reduction methods. The most classical result is a bound for balanced truncation and singular perturbation approximation that is applicable for asymptotically stable linear time-invariant systems with homogeneous initial conditions. Recently, there have been a few attempts to generalize the balancing-related reduction methods to the case with inhomogeneous initial conditions, but the existing error bounds for these generalizations are quite restrictive. Particularly, it is required to restrict the  initial conditions to a low-dimensional subspace, which has to be chosen before the reduced model is constructed. In this paper, we propose an estimator that circumvents this hard constraint completely.  Our estimator is applicable to a large class of reduction methods, whereas the former results were only derived for certain specific methods. Moreover, our approach yields to significantly more effective error estimation, as also will be demonstrated numerically.
\end{abstract}
\end{frontmatter}

\section{Introduction}
We consider linear time-invariant systems with inhomogeneous initial conditions
\begin{align}
\begin{aligned}\label{eq:fom}
 \dot{\bv x} (t) &= \bt A \bv x(t) + \bt B \bv u(t)\\
   \bv y(t) &= \bt C \bv x(t), \hspace{0.55cm} \bv x(0) = \bv x_0,
\end{aligned}  
\end{align}
with  $\bt A \in \R^{N,N}$, $\bt B \in \R^{N,q}$ and $\bt C \in \mathbb{R}^{N,p}$. The matrix $\bt A$ is assumed to be Hurwitz, 
which is equivalent to the asymptotic stability of the system. We refer to \eqref{eq:fom} as the full order model (\acro{FOM}) and assume that the state $\bv x$ is high-dimensional compared to the dimension of the input $\bv u$ and output $\bv y$ ($N\gg p,q$). The latter motivates the use of model reduction, which is a methodology to construct a surrogate model with lower computational costs, the so-called reduced model (\acro{ROM}). The standard projection-based ansatz relies on Petrov-Galerkin projection, which is realized with appropriate reduction bases $\bt V, \bt W \in \mathbb{R}^{N,n}$, $n \ll N$. Under the usual bi-orthogonality  assumption $\bt W^T \bt V = \bt I_n$, the \acro{ROM} then takes the form
\begin{align}
\begin{aligned}\label{eq:rom}
 \dot{\rom{\bv{x}}} (t) &= \rom{\bt{A}} \rom{\bv{x}}(t)+ \rom{\bt B} \bv u(t)\\
   \rom{\bv y}(t) &= \rom{\bt C} \rom{\bv x}(t),\hspace{0.55cm} \rom{\bv x}(0) = \bt W^T \bv x_0,
\end{aligned}
\end{align}
with $\rom{\bt{A}} = \bt W^T \bt A \bt V$, $\rom{\bt{B}} = \bt W^T \bt B$ and $\rom{\bt{C}} = \bt C \bt V$. The foremost aim is a high fidelity approximation of the output, $\bv y(t) \approx \rom{\bv{y}}(t)$, for all scenarios of interest.

Balancing-related model reduction methods, such as balanced truncation (\acro{BT}) or singular perturbation approximation, are powerful model reduction methods. Originally, they were derived for linear systems with homogeneous initial conditions, i.e., \eqref{eq:fom} with $\bv x_0 = \bv 0$. In this setting, a priori error bounds can be shown \cite{morMoo81,book:antoulas2005,book:LinearRobustControl}. But the treatment of inhomogeneous initial conditions is non-standard for balancing-related and other system-theoretic methods, which represents a limitation for their practical use in the time-domain analysis. Recently, a few attempts have been made to alleviate this limitation and to implement inhomogeneous initial conditions in system-theoretic methods \cite{art:Daraghmeh2019,art:morHeiRA11,morBeaGM17,art:bt-schroeder20}. The basic idea these works pursue is a training of the models towards a predefined space of expected initial conditions. Error bounds have been derived for the respective methods, but they strictly require the initial conditions to lie in the space the reduced models have been trained for.

In this paper, we propose an error estimator that is applicable for arbitrary initial conditions, i.e. completely circumvents the necessity of posing restrictions or assuming prior knowledge on the simulation setups. Consequently, our estimator can be used to certify or reject a \acro{ROM} on the fly for any given scenario, without the need of costly recomputations. Moreover, as we will demonstrate, our estimator is more effective than the existing ones, in the sense that it estimates the error significantly sharper. In addition, it is applicable to any asymptotically stable \acro{ROM}, whereas the ones from \cite{art:Daraghmeh2019,art:morHeiRA11,morBeaGM17,art:bt-schroeder20} can only be used in combination with the specific method considered there. The results in the cited works strongly rely on a control-type viewpoint and the reinterpretation of the initial condition as an additional impulsive input. We consider the error estimation problem as an observability problem instead. More specifically, the reduction error is split into a controlled part (independent of the initial conditions) and an uncontrolled part (independent of the inputs) and derive an explicit representation of the uncontrolled reduction error in terms of the observability Gramian. This yields a perfect error prediction when the Gramian can be calculated exactly and no other round-off errors are present. When only an approximation on the Gramian is available--certainly the more realistic scenario in a large-scale setting--our approach yields an error estimator that is as effective as the Gramian approximation omits. 
Our approach can be implemented in an online-efficient way, without any restrictions on the initial conditions. The required offline computations are also moderate given the (approximate) observability Gramian is available.
This makes our approach particularly well-suited for the combination with the augmented \acro{BT} approach \cite{art:morHeiRA11} as well as for the splitted reduction approach \cite{morBeaGM17}, in which the controlled and uncontrolled parts are separately reduced. For the splitted approach, we assume that \acro{BT} is used for the controlled part, which implies that the Gramian has to be determined anyway. But, notably, no restrictions on the reduction of the uncontrolled part have to be posed apart from asymptotic stability of the resulting \acro{ROM}.

The remainder of the paper is structured as follows. The system-theoretic concepts needed in this paper are outlined in Section~\ref{sec:system-theory}. The \acro{BT} method for systems with homogeneous initial conditions and the use of low-rank approximation in this context are briefly explained in Section~\ref{sec-bt-controlled}.  The main result of this paper, i.e., our proposed error estimator related to inhomogeneous initial conditions, is presented in Section~\ref{sec:error-bound}. Then different reduction methods that account for inhomogeneous initial conditions are recapitulated from literature (Section~\ref{sec:mor-inhom}). The effectiveness and wide applicability of our estimator is demonstrated in Section~\ref{sec:num-results} using different methods and examples, particularly also a larger-scale problem where the required Gramian is only approximately determined.

\section{System-theoretic concepts} \label{sec:system-theory}
Let the input $\bv u$ be square-integrable in time $t\in [0,\infty)$. Then the respective output of the \acro{FOM} \eqref{eq:fom} is also square-integrable, and its $\funSpace{L}^2$-norm is defined by
\begin{align*}
	||\bv  y||_{\funSpace{L}^2} =  \sqrt{\int_{0}^\infty ||\bv y(t)||^2 dt }.
\end{align*}
Here and in the following, $ ||\cdot ||$ denotes the Euclidean vector norm. The few system-theoretic results we utilize in this paper are summarized in the upcoming, see \cite{book:antoulas2005,book:LinearRobustControl} for details.
As the system matrix $\bt A$ of \eqref{eq:fom} is assumed to be Hurwitz, the observability Gramian $\Ob \in \R^{N,N}$ and the controllability Gramian $\bt{P}\in \R^{N,N}$ are well-defined, symmetric positive semi-definite  and given as the unique solutions of the two Lyapunov equations
\begin{align} \label{eq:lyap}
	\bt A^T \Ob + \Ob \bt A =- \bt C^T \bt C, \hspace{0.3cm} \bt A \bt{P} + \bt P \bt A^T =- \bt B \bt B^T.
\end{align}
In particular, the observability Gramian has the representation
\begin{align*}
 \Ob = \int_{0}^\infty \mexp{t \bt A^T} \bt C^T \bt C \mexp{t \bt A} dt.
\end{align*}
The output energy for a system of form \eqref{eq:fom} with trivial input, $\bv u(t) = 0$ for $t\geq 0$, is given by
\begin{align}\label{eq:out-energy}
		 	||\bv  y||_{\funSpace{L}^2}^2  =  \int_{0}^\infty  (\bt C \mexp{t \bt A} \bv x_0)^T (\bt C \mexp{t \bt A} \bv x_0) dt = \bv x_0^T  \Ob \bv x_0.
\end{align}
Thus, $\Ob$ describes a measure for the observability of a given (initial) state $\bv x_0$. Similar considerations motivate that $\bt P$ induces a measure for controllability. The Gramians depend on the given state representation and, consequently, change under state transformations. But the eigenvalues of the matrix $(\bt P \Ob)$ are invariants of the system. Their square-roots $\sigma_1 \geq \ldots \geq \sigma_N \geq 0$, known as the Hankel singular values, play a crucial role in system-theoretic model reduction.

The input-output map $\bv u \mapsto \bv y$ of a system with homogeneous initial conditions (i.e., \eqref{eq:fom} with $\bv x_0 = \bv 0$) is also invariant under state transformations. In frequency space it is characterized by the transfer function, $\bt G(s) = \bt C (s\bt I - \bt A)^{-1} \bt B$ for $s\in \mathbb{C}$. According to the Plancherel theorem, the $\funSpace{L}^2$-norms in time and frequency space are equal, which implies the bound
\begin{align*}
|| \bv y||_{\funSpace{L}^2} \leq || \bt G||_{\funSpace{H}_\infty} || \bv u||_{\funSpace{L}^2}.
\end{align*}
The expression $|| \bt G||_{\funSpace{H}_\infty}$ refers to the $\funSpace{H}_\infty$-norm of the transfer function, which is the operator-norm induced by the $\funSpace{L}^2$-norm for the input-output map.

\section{(Approximate) balanced truncation} \label{sec-bt-controlled}

{\acro{BT}} is a widely used system-theoretic model reduction method for systems with homogeneous initial conditions. It provably preserves asymptotic stability and, even more important in practice, comes with an a priori error bound. Here, we summarize the rudimentary ideas of the standard method \cite{morMoo81} and its approximate version \cite{book:dimred2003}, which is typically used in the large-scale setting.



\subsection{Balanced truncation} \label{subsec:btst}
The basic idea behind \acro{BT} is to construct reduced models that relate to a truncation of the \acro{FOM} in a balanced form. Balanced in this context refers to a realization, for which the observability and controllability Gramians are both equal to the diagonal matrix that has the Hankel singular values $\sigma_1 \geq \ldots \geq \sigma_N \geq 0$ in descending order as diagonal entries. For ease of presentation, let us assume $\sigma_n > \sigma_{n+1}$. Then the output $\rom{\bv y }$ of the \acro{ROM} of order $n$ that is obtained by {\acro{BT}} can be shown to fulfill
\begin{align}\label{eq:bthom-bound}
		|| \bv y - \rom{\bv y } ||_{\funSpace{L}^2} &\leq {\alpha} || \bv u||_{\funSpace{L}^2}, \qquad  \text{with } \, {\alpha} = 2\sum_{i=n+1}^N \sigma_i
\end{align}
for any square-integrable input $\bv u$. This bound can also be specified as $|| \bt G -\rom{\bt G}||_{\funSpace{H}_\infty}\leq {\alpha}$, whereby $\rom{\bt G}$ denotes the transfer function of the \acro{ROM}, see \cite{book:antoulas2005,morMoo81} for details.

\subsection{Low-rank approximation of Gramians} \label{subsec:low-rank}
In order to perform the full balancing, the Gramians $\Ob$ and $\bt P$ need to be determined. This is by far the most computational demanding step in the implementation of {\acro{BT}}. When the system matrix $\bt A$ is high-dimensional and sparse, a typically much more efficient and numerically robust approach is obtained by replacing the exact Gramians with low-rank approximations, e.g., for the observability Gramian
\begin{align} \label{eq:Ob-lowrank}
	\Ob \approx \bt U^T \bt U, \quad \bt U \in \R^{m,N}, \quad m\ll N.
\end{align}
The factors $\bt U$ are constructed such that $\bt U^T \bt U$ solve the Lyapunov equations \eqref{eq:lyap} up to a small approximation error. Note that the product never needs to be formed explicitly in the implementation, but all operations can be directly performed on $\bt U$. The ADI-based and projection-based approaches \cite{art:Simoncini2016,SaaKB-mmess} have shown to be very reliable and efficient in determining low-rank factors that yield high fidelity approximations. As also shown in the context of model reduction \cite{book:dimred2003,morAntBG10}, the low-rank errors are negligible in many practical examples. But it should be mentioned that the asymptotic stability preserving property and the bound \eqref{eq:bthom-bound} cannot be proven by the standard results on \acro{BT} any longer when not the exact Gramians are used. The error analysis revolving around low-rank solvers is a topic on its own, cf. \cite{art:morAntSZ02,art:Simoncini2016}, and not further addressed in this paper.

\section{Proposed error estimator for inhomogeneous initial conditions} \label{sec:error-bound}

This section presents our main result. We show that the reduction error related to the initial conditions can be effectively estimated for a large class of reduced models, given the (approximate) observability Gramian is available. 
As a preliminary step, a splitting of the {FOM}, respectively the reduction error, into a control-dependent part and an uncontrolled part is needed (Section~\ref{subsec:split-romerror}). The estimator itself is derived in Section~\ref{subsec:estimator-derivation}.

\subsection{Splitting of the reduction error} \label{subsec:split-romerror}
Starting point for our considerations is the following well-known splitting of the \acro{FOM} output, cf. \cite{morBeaGM17,book:dimred2003},
\begin{subequations} \label{eq:superpos}
\begin{align}
\begin{aligned}
	\bv y(t) &= \underbrace{\bt C e^{t \bt A}\bv x_0 }_{\bv y_{x_0}(t)}  + \underbrace{\int_{0}^t \bv C e^{(t-\tau) \bt A} \bt B \bv u(\tau) d\tau}_{\bv y_u(t)}\\
	 & \bv y_{x_0}(t) = \bv C \bv z(t), \hspace{0.3cm} \bv y_{u}(t) = \bv C \bv x_u(t).
\end{aligned}
\end{align}
The related split state representations read
\begin{align}
		\dot{\bv{z}}(t) &= \bt A \bv z(t), \hspace{0.6cm}  &  \bv z(0) = \bv x_0,  \label{eq:superpos-uc} \\
		\dot{\bv{x}}_u(t) &= \bt A \bv x_u(t) + \bt B \bv u(t), \hspace{0.3cm} & \bv x_u(0)= \bv 0  \label{eq:superpos-hom}
\end{align}
\end{subequations}
Particularly, \eqref{eq:superpos-hom} describes the controlled part \textit{(with homogeneous initial conditions)}, and \eqref{eq:superpos-uc} the uncontrolled part \textit{(with inhomogeneous initial conditions)}.

We assume that the reduced output $\rom{\bv y}$ allows for a similar split representation,
\begin{subequations} \label{eq:superpos-r}
\begin{align}
		\rom{\bv y}(t) &= \rom{\bv y}_{x_0}(t) + \rom{\bv y}_{u}(t)  =  \rom{\bv C} \rom{\bv z}(t) + \rom{\bt C}_u \rom{\bv x}_u(t)
\end{align}
with the uncontrolled and a controlled reduced part described by
\begin{align}
		\dot{\rom{\bv z}}(t) &= \rom{\bt A} \rom{\bv z}(t), \hspace{0.2cm}  &\rom{\bv z}(0) = \bt W^T \bv x_0,  \label{eq:superpos-r-uc} \\
		\dot{\rom{\bv x }}_u(t) &= \rom{\bt A}_u \rom{\bv x}_u(t) + \rom{\bt B}_u \bv u(t), \hspace{0.3cm}  & \rom{\bv x}_u(0)= \bv 0. \label{eq:superpos-r-hom}
\end{align}
\end{subequations}
For the standard approach according to \eqref{eq:rom}, the reduction bases $\bt W$, $\bt V$ are used to reduce the system as a whole, or equivalently \eqref{eq:superpos-r-uc}-\eqref{eq:superpos-r-hom} are both constructed with the same bases (i.e., $\rom{\bt A}= \rom{\bt A}_u = \bt W^T \bt A \bt V$, $\rom{\bt B}_u = \bt W^T \bt B$ and $\rom{\bt C} = \rom{\bt C}_u = \bt C \bt V$). Let us stress that our results are also valid in the more general setting that the reduction of the controlled and the uncontrolled is performed separately (i.e., $\rom{\bt A}\neq \rom{\bt A}_u$ and $\rom{\bt C} \neq \rom{\bt C}_u$), even when certain non-projection-based approaches are used,  e.g., the singular perturbation approximation.

Subtracting the split \acro{FOM} and \acro{ROM} from each other and applying the triangle inequality yields the estimate
\begin{align} \label{eq:error-splitting}
	||\bv y - \rom{\bv y}||_{\funSpace{L}^2} \leq ||\bv y_{x_0} -\rom{\bv y}_{x_0}||_{\funSpace{L}^2} + ||\bv y_{u} -\rom{\bv y}_{u}||_{\funSpace{L}^2}
\end{align}
for the reduction error. The term $||\bv y_{u} -\rom{\bv y}_{u}||_{\funSpace{L}^2}$ is the one that is well-studied for system-theoretic model reduction. Given that a balancing-related method is used for it, which we assume throughout the paper, the a priori error bounds from literature can be applied. Thus, we are only concerned with $||\bv y_{x_0} -\rom{\bv y}_{x_0}||_{\funSpace{L}^2}$, the reduction error of the uncontrolled part. The error estimator we derive is well-posed under the following assumption.
\begin{assum} \label{assum:fomromHurwitz}
The matrices $\bt A$ and $\rom{\bt A}$ of the \acro{FOM} \eqref{eq:superpos-uc} and the \acro{ROM} \eqref{eq:superpos-r-uc} are both assumed to be Hurwitz.
\end{assum}
The estimator depends explicitly on the (approximate) observability Gramian of the \acro{FOM}. This further motivates to combine our approach with a balancing-related method, applied to either reduce the full system or to at least one of its parts. In this situation, the additional computational costs for the estimator are moderate, as the Gramian has already been determined during the construction of the \acro{ROM}.

\subsection{Derivation of proposed estimator} \label{subsec:estimator-derivation}
We derive an estimator for $||\bv y_{x_0} -\rom{\bv y}_{x_0}||_{\funSpace{L}^2}$, the reduction error related to the initial conditions.
Subtracting the uncontrolled parts of the \acro{FOM}  and the \acro{ROM}, i.e. \eqref{eq:superpos-uc}  and \eqref{eq:superpos-r-uc}, we obtain
\begin{align}
\begin{aligned} \label{eq:errfactor}
\dot{\bv{z}}(t) &= \bt A \bv z(t) , \hspace{0.5cm} \bv z(0) = \bv x_0\\
\dot{\rom{\bv z}}(t) &= \rom{\bt A} \rom{\bv z}(t), \hspace{0.5cm} \rom{\bv z}(0) = \bt W^T \bv x_0\\
 \bv y_{x_0}(t)-\rom{\bv y}_{x_0}(t)&= \bv C \bv z(t) -\rom{\bv C} \rom{\bv z}(t).
\end{aligned}
\end{align}
Equation \eqref{eq:errfactor} describes an uncontrolled system in the extended state
$\cvv{z} = [\bv z^T, \rom{\bv{z}}^T ]^T$, which can also be stated as
\begin{align*}
\begin{aligned} 
 \dot{\cvv z}(t) &= \ct A \cvv z(t), \hspace{0.5cm} \cvv z(0) = \cvv z_0,
 \hspace{0.7cm}
\bv y(t)- \rom{\bv y}(t)= \ct C \cvv z(t) 
\\ 
 \ct{A}& =
 \begin{bmatrix}
 	\bt A &  \\
 			& \rom{\bt{A}}
 \end{bmatrix},  
 \hspace{0.4cm}
 \cvv z_0 = \begin{bmatrix}
 \bv x_0 \\
  \bt W^T \bv{x}_0
\end{bmatrix},
\hspace{0.5cm}
\ct C = \begin{bmatrix}
 \bt C, & \bt  -\rom{\bt C} 
 \end{bmatrix}.   
\end{aligned}
\end{align*}

Under Assumption~\ref{assum:fomromHurwitz}, the observability Gramian $\Obc$ of the error system is well-defined, symmetric positive semi-definite and the solution of the Lyapunov equation
\begin{align*}
	\ct A^T \Obc + \Obc \ct A = - \ct C^T \ct C,   \hspace{0.5cm} 
	\Obc = \begin{bmatrix}
	\Ob &  \mmark{\Ob} \\
		\mmark{\Ob}^T	&  \rom{\Ob}
	\end{bmatrix} \in \R^{N+n,N+n}.
\end{align*}
The latter is equivalent to the three equations
\begin{subequations}\label{eq:lya-eb} 
\begin{align}	
	\bt A^T \Ob + \Ob \bt A &=- \bt C^T \bt C  \label{eq:lya-eb-a} \\
	\bt A^T \mmark{\Ob} + \mmark{\Ob} \rom{\bt A} &= \hspace{0.3cm} \bt C^T \rom{\bt C} \label{eq:lya-eb-b} \\
	\rom{\bt A}^T \rom{\Ob} + \rom{\Ob} \rom{\bt A} &=-\rom{\bt C}^T \rom{\bt C}. \label{eq:lya-eb-c}
\end{align}
\end{subequations}
Applying the notion of observability to \eqref{eq:errfactor} and defining $\rom{\bv x}_0 = \bt W^T \bv x_0$ yields
\begin{align}
\begin{aligned}\label{error-gen}
	||\bv y_{x_0} - \rom{\bv y}_{x_0}||_{\funSpace{L}^2}^2 &= \cvv z_0^T \Obc \cvv z_0 \\
	&= \bv x_0^T \Ob \bv x_0 + 2 \bv x_0^T \mmark{\Ob} \rom{\bv x}_0 + \rom{\bv x}_0^T \rom{\Ob} \rom{\bv x}_0,
\end{aligned}
\end{align}
cf. Section \ref{sec:system-theory}. By replacing the observability Gramian $\Ob$ of the \acro{FOM} with a low-rank approximation, we obtain the following error bound and error estimator, respectively.

\begin{thm}[Error bound and estimator]\label{thrm:errboundxo}
Let Assumption~\ref{assum:fomromHurwitz} hold. Let $\Ob$, $\mmark{\Ob}$ and $\rom{\Ob}$ be the solutions to \eqref{eq:lya-eb}, and let $\bt U \in \R^{m,N}$ with $m\leq N$ define the approximation $\Ob \approx \bt U^T \bt U$. Then the output error defined by \eqref{eq:errfactor} fulfills
	\begin{align*}
	||\bv y_{x_0} - \rom{\bv y}_{x_0}||_{\funSpace{L}^2} &= \sqrt{\Delta_{\bv x_0}^2 + \bv x_0^T (\Ob-\bt U^T \bt U)\bv x_0}
	\\ &\leq \sqrt{\Delta_{\bv x_0}^2 + ||\Ob-\bt U^T \bt U||_2  ||\bv x_0||^2}, \hspace{0.6cm} 
	\end{align*}
whereby $\Delta_{\bv x_0}$ reads
\begin{align*}
\Delta_{\bv x_0}= \sqrt{  ||\bt U \bv x_0||^2 + 2 \bv x_0^T \mmark{\Ob} \bt W^T {\bv x}_0 + \bv x_0^T \bt W \rom{\Ob}\bt W^T  \bv x_0 }.
\end{align*}
\end{thm}
\begin{pf}
It holds $\bv x_0^T \Ob \bv x_0 = \bv x_0^T \Ob_a \bv x_0 + \bv x_0^T (\Ob -\Ob_a) \bv x_0$ for any $\Ob_a \in \R^{N,N}$. Using this for $\Ob_a= \bt U^T \bt U$ and using \eqref{error-gen} the assertion can be derived.
\end{pf}
We propose $\Delta_{\bv x_0}\approx ||\bv y_{x_0} - \rom{\bv y}_{x_0}||_{\funSpace{L}^2}$ as error estimation. Our estimator is exact (up to small round-off errors) when $\bt U$ is chosen as the Cholesky factor of $\Ob$, which can be done in a small-scale setting. In other cases, the effectivity of our estimator is determined by the fidelity of the Gramian approximation. Fortunately, the low-rank solvers yield very high-fidelity approximations in most relevant settings, cf. Section~\ref{subsec:low-rank}.
Our estimator allows for an efficient offline-online decomposition, by simply precalculating the matrices $\bt U$,  $\mmark{\Ob}$ and $\rom{\Ob}$. As we assume a balancing-related method to be used (cf. Section~\ref{sec:mor-inhom}),  $\bt U$ is already determined in the reduction step. The calculation of the other two matrices is less computationally demanding. Particularly, $\rom{\Ob}$ is the solution of a low-dimensional Lyapunov equation and can be obtained without difficulty. To obtain $\mmark{\Ob}$, the sparse-dense Sylvester equation \eqref{eq:lya-eb-b} is solved using the algorithm proposed in \cite{inproc:Vasilyev05,inproc:PoussotVassal2011}: First the complex Schur decomposition of $\rom{\bt A}$ is determined, and then a small number of equation solves of dimension $N$ are performed. 
\begin{rem}[Relation to  $\funSpace{H}_2$-optimization] \label{rem:sparse-dense-sylvester}
The equations \eqref{eq:lya-eb} also play a crucial role in the iterative $\funSpace{H}_2$-optimal model reduction method \acro{IRKA} \cite{morAntBG10} and its so-called TSIA implementation \cite{art:TSIA-Xu11}. The iterative procedures these algorithms consist of implicitly require solving equations of the form \eqref{eq:lya-eb-b}-\eqref{eq:lya-eb-c} multiple times.
\end{rem}

\section{System-theoretic model reduction with inhomogeneous initial conditions} \label{sec:mor-inhom}

The existing extensions of system-theoretic reduction methods for inhomogeneous initial conditions \cite{art:morHeiRA11,art:Daraghmeh2019,morBeaGM17,art:bt-schroeder20} rely on a training towards a user-defined set of expected initial values. A training matrix $\bt X_0 \in \R^{N,N_0}$ with $N_0\ll N$ is chosen, and the expected initial values are those that can be written as
\begin{align} \label{eq:x0restriction}
\bv{x}_0 = \bt X_0 \bv v_0, \qquad \text{for } \bv v_0 \in \R^{N_0}.
\end{align}
The impact of $\bv x_0$ on the dynamics is then formally rewritten as an additional impulsive input, so that a system-theoretic methods can be used. We visit two different approaches that follow this scheme and comment on the error bounds proposed in former works.

\subsection{Reduction relying on augmentation} \label{subsec-BTaug}
The method from \cite{art:morHeiRA11} is abbreviated as {\acro{BT-aug}} in the upcoming. It is based on the augmented system
\begin{align*}
	\dot{\bv{x}}(t) &= \bt A \bv x(t) + \bt B_{\text{aug}} \bv u_{\text{aug}}(t)\\
	\bv y_{\text{aug}}(t)&= \bt C \bv x(t), \hspace{0.55cm}\bv x(0) =\bv 0,
\end{align*}
 with $\bt B_{\text{aug}} = [\bt B, \bt X_0]$, which has homogeneous initial conditions and is formally equivalent to the \acro{FOM} \eqref{eq:fom} when $\bv u_{\text{aug}}(t)= [\bv u(t)^T,\bv v_0^T\delta(t)]^T$ is chosen, with $\delta$ denoting the Dirac impulse. Reduction bases $\bt V$, $\bt W$ are constructed by using \acro{BT} on the augmented system and the \acro{ROM} is obtained following the standard projection ansatz \eqref{eq:rom}. It should be mentioned that $\bv u_{\text{aug}}$ is not square-integrable, and thus the standard error bound for \acro{BT} (i.e., \eqref{eq:bthom-bound}) is not applicable. The reference \cite{art:morHeiRA11} derives an alternative error bound under the restriction \eqref{eq:x0restriction} on the initial value $\bv x_0$. It is of the form
 	\begin{align}\label{eq:bound-BTaug}
		||\bv y- \rom{\bv y}_{\text{\acro{BT-aug}}}||_{\funSpace{L}^2} &\leq \alpha_{\mathrm{aug}} ||\bv u||_{\funSpace{L}^2} + \beta_{\mathrm{aug}} || \bv v_0 ||,
	\end{align}
whereby $\alpha_{\mathrm{aug}}$ depends on the truncated Hankel singular values for the augmented system. The constant  $\beta_{\mathrm{aug}}$ depends nonlinearly on $\alpha_{\mathrm{aug}}$ and the matrices $\bt X_0$ and $\bt A$, we refer to \cite{art:morHeiRA11} for details.

\subsection{Reduction relying on splitting} \label{subsec-BTBT}
The approach suggested in \cite{morBeaGM17} is based on the splitting \eqref{eq:superpos} of the \acro{FOM}. It reduces the controlled part and the uncontrolled part individually. The use of \acro{BT} is proposed for the controlled part \eqref{eq:superpos-hom}. For the reduction of the uncontrolled part \eqref{eq:superpos-uc}, the auxiliary system
\begin{align*}
	\dot{\bv{z}}(t) &= \bt A \bv z(t) + \bt B_{x_0} \bv u_{x_0}(t)\\
	\bv y_{x_0}(t) &= \bt C \bv z(t), \hspace{0.55cm}  \bv z(0) = \bv 0
\end{align*}
with $\bt B_{x_0} = \bt X_0$ is considered, which is formally equivalent to the uncontrolled system using $\bv u_{x_0}(t)= \bv v_0 \delta(t)$, $t\geq0$. The reference \cite{morBeaGM17} suggests to reduce this part by either {\acro{BT}} or \acro{IRKA}. The method variant that uses \acro{BT} for the uncontrolled part is abbreviated as \acro{BT-BT} from here on. It guarantees asymptotic stability, and an error bound for it is stated in the reference, which holds under restriction \eqref{eq:x0restriction} on $\bv x_0$. It takes the form 
 	\begin{align}\label{eq:bound-BTBT}
				||\bv y- \rom{\bv y}_{\text{\acro{BT-BT}}}||_{\funSpace{L}^2} &\leq \alpha_{\mathrm{spl}} ||\bv u||_{\funSpace{L}^2} + \beta_{\mathrm{spl}} || \bv v_0 ||,
	\end{align}
with $\alpha_{\mathrm{spl}}$ obtained by the standard result on \acro{BT} (Section~\ref{subsec:btst}). Note that the explicit determination of $\beta_{\mathrm{spl}}$ requires a fully balanced realization of the uncontrolled \acro{FOM} part \eqref{eq:superpos-uc}, which  makes it only feasible in a small-scale setting.

The other method variant using \acro{IRKA} for the uncontrolled part has the disadvantage that no a priori stability guarantees can be given. We therefore propose an asymptotic stability preserving modification, namely to use \acro{ISRK} \cite{art:isrk-Gugercin08,inproc:PoussotVassal2011}  for the reduction of the uncontrolled part. This algorithm is very similar to \acro{IRKA}. The difference lies in the additional constraint $\bt W = \Ob \bt V (\bt V^T \Ob \bt V)^{-1}$, which is used in \acro{ISRK} to determine $\bt W$ from $\bt V$ in every step of the underlying iteration. Note that \acro{ISRK} is computationally more demanding than \acro{IRKA} when used as a standalone reduction method, because of the necessity of the Gramian  $\Ob$. However, for the splitted reduction approach considered here, the latter is no concern, as the Gramian has already been determined to perform \acro{BT} for the controlled part. 


\subsection{Comparison of proposed and existing error bounds}

In case the Gramians can be determined exactly, our estimator yields a strict bound, which is applicable under the same assumptions as the bounds for \acro{BT-aug}  \eqref{eq:bound-BTaug} and \acro{BT-BT} \eqref{eq:bound-BTBT} from literature. But our result is also applicable to other \acro{ROM}s and under more general settings. Particularly, it does not require the restriction \eqref{eq:x0restriction} on the initial value $\bv x_0$.
\begin{cor} \label{cor:strict-bound}
	Let the assumptions of Theorem~\ref{thrm:errboundxo} hold, and let additionally $\Ob = \bt U^T \bt U$, i.e. $\bt U$ be the exact Cholesky factor of $\Ob$. Then the output error \eqref{eq:error-splitting} for a \acro{ROM} constructed by \acro{BT-aug} or by the splitted method (\acro{BT-BT} or another method from Section~\ref{subsec-BTBT}) fulfills the error bound
	\begin{align*}
		||\bv y- \rom{\bv y}||_{\funSpace{L}^2} \leq \tilde{\alpha} ||\bv u||_{\funSpace{L}^2} + \Delta_{\bv x_0},
	\end{align*}
	with $\tilde{\alpha} =\alpha_{\mathrm{aug}}$ for \acro{BT-aug}, and $\tilde{\alpha} =\alpha_{\mathrm{spl}}$ for splitted \acro{ROM}s.
\end{cor}

For a more illustrative comparison of our estimator and the ones from literature, we state the following corollary.
\begin{cor}
	Let the assumptions of Theorem~\ref{thrm:errboundxo} hold, and let $\bar{\bv X}_0\in \R^{n,\bar{n}_0}$ and $\bv v_0 \in \R^{\bar{n}_0}$. Then our error estimator fulfills for $\bv x_0 = \bar{\bv X}_0 \bv v_0$ that
	\begin{align*}
		\Delta_{\bv x_0} =  \sqrt{ \bv v_0^T \bt Z \bv v_0} \leq \sqrt{|| \bt Z ||_2} \, ||\bv v_0||,
	\end{align*}
	where $\bt Z \in \R^{\bar{n}_0,\bar{n}_0}$ is given by
		\begin{align*}
	\bt Z =&  \bar{\bt X}_0^T \left( \bt U^T \bt U +  \mmark{\Ob} \bt W^T +  \bt W \mmark{\Ob}^T  +  \bt W\rom{\Ob}\bt W^T    \right) \bar{\bt X}_0.
	\end{align*}	

\end{cor}
The corollary follows straight forwardly from Theorem~\ref{thrm:errboundxo}. Notably, our result is sharp in the sense that it holds $\Delta_{\bv x_0} = \sqrt{|| \bt Z ||_2} \, ||\bv v_0||$ for $\bv x_0 = \bar{\bv X}_0 \bv v_0$ with $\bv v_0$ chosen as an eigenvector of $\bt Z$ to the biggest eigenvalue. 
Let us also stress again that our estimator is applicable for arbitrary initial values and an efficient offline-online decomposition that is independent of the initial values is possible, based on the formula stated in Theorem~\ref{thrm:errboundxo}.

\section{Numerical results} \label{sec:num-results}

Our numerical studies aim to illustrate the effectivity of our error estimator, which holds independently of the reduction method as well as in the small- and the large-scale setting. 
In Section~\ref{subsec:num-beam}, a small-scale problem from literature is used and direct comparisons to the error estimators from literature are made. The efficiency and reliability in the presence of low-rank approximation errors in the observability Gramian is demonstrated at a larger-scale problem in Section~\ref{subsec:num-convdiff}.

All numerical results have been generated using \texttt{MATLAB} Version 9.1.0 (R2016b) on an Intel Core i5-7500 CPU with 16.0GB RAM and the \texttt{MESS} toolbox \cite{SaaKB-mmess} for the sparse solvers and the reduction routines, which are partly based on \cite{code:voigt20-btinhom}.
For better reproducibility, the code and the benchmark data are provided in \cite{code:bls22-btinhom}.


The reduction errors up to a time $t>0$ are defined as
\begin{align*}
	E(t) &= \sqrt{\int_0^t ||\bv y(\tau)- \rom{\bv y}(\tau)||^2 d\tau}\\
		E_{x_0}(t) &= \sqrt{\int_0^t ||\bv y_{x_0}(\tau)- \rom{\bv y}_{x_0}(\tau)||^2  d\tau}
\end{align*}
for the output $\bv y$ and its uncontrolled part $\bv y_{x_0}$, respectively. By construction, $E(t) \leq  ||\bv y- \rom{\bv y}||_{\funSpace{L}^2}$ for $t\geq 0$, as the latter is the limit for $t\rightarrow \infty$. For uncontrolled scenarios, it holds $E(t)= E_{x_0}(t)$. We evaluate the errors by simulating the models with the \texttt{MATLAB} built-in integrator \texttt{ode45} 
and then approximating the integrals with a trapezodial rule. Given a final simulation time $T$, these numerical approximations are performed on a time mesh that is highly resolved around the origin, consisting of $10\,000$ logarithmically distributed points between $10^{-20}T$ and $T$, and zero as starting point.

\subsection{Beam equation} \label{subsec:num-beam}
We consider the beam example from the SLICOT benchmark collection \cite[Section~24]{book:dimred2003}, which is also used in \cite{art:bt-schroeder20}.  The system matrices have the dimensions $n=349$ and $q=p=1$, i.e., it is a small-scale problem with a single input and a single output. As in \cite{art:bt-schroeder20}, we choose $\bt X_0 = [X_0^{(i,j)}] \in \R^{n,2}$ with $X_0^{(5,1)}=1$, $X_0^{(101,2)}=100$ and zeros elsewhere as training matrix for the \acro{ROM}s. For the simulation, we consider two different setups.
\begin{itemize}
\item \textbf{{Trained case:}} The initial value $\bv x_0= \bt X_0 \bv v_0$ is chosen with $\bv v_0 = [10,-1]^T$. The final simulation time is given by $T = 1000$, and the input by
\begin{align*}
	u(t) = 
	\begin{cases}
	1, & t \in [100,200] \\
	0, & \text{else.}
	\end{cases}
\end{align*}
\item \textbf{{Not trained case:}} The initial condition and end time are set to $\bv x_0 = [5,5,\ldots,5]^T $ and $T = 10\, 000$, respectively. No input is used, yielding an uncontrolled scenario.
\end{itemize}

\begin{figure}[tb]
\begin{center}
\includegraphics[height=3.1cm]{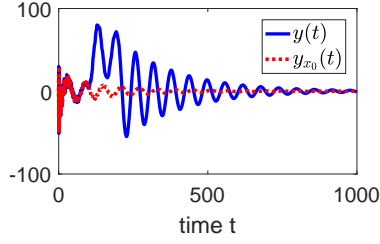}    
\caption{Beam, trained case. \acro{FOM} output $y(t)$ and uncontrolled part $y_{x_0}(t)$.  
\label{fig:beam-trainedsol}
}                                                           
\end{center}                            
\end{figure}

\begin{figure}[tb]
\begin{center}
\includegraphics[height=4.1cm]{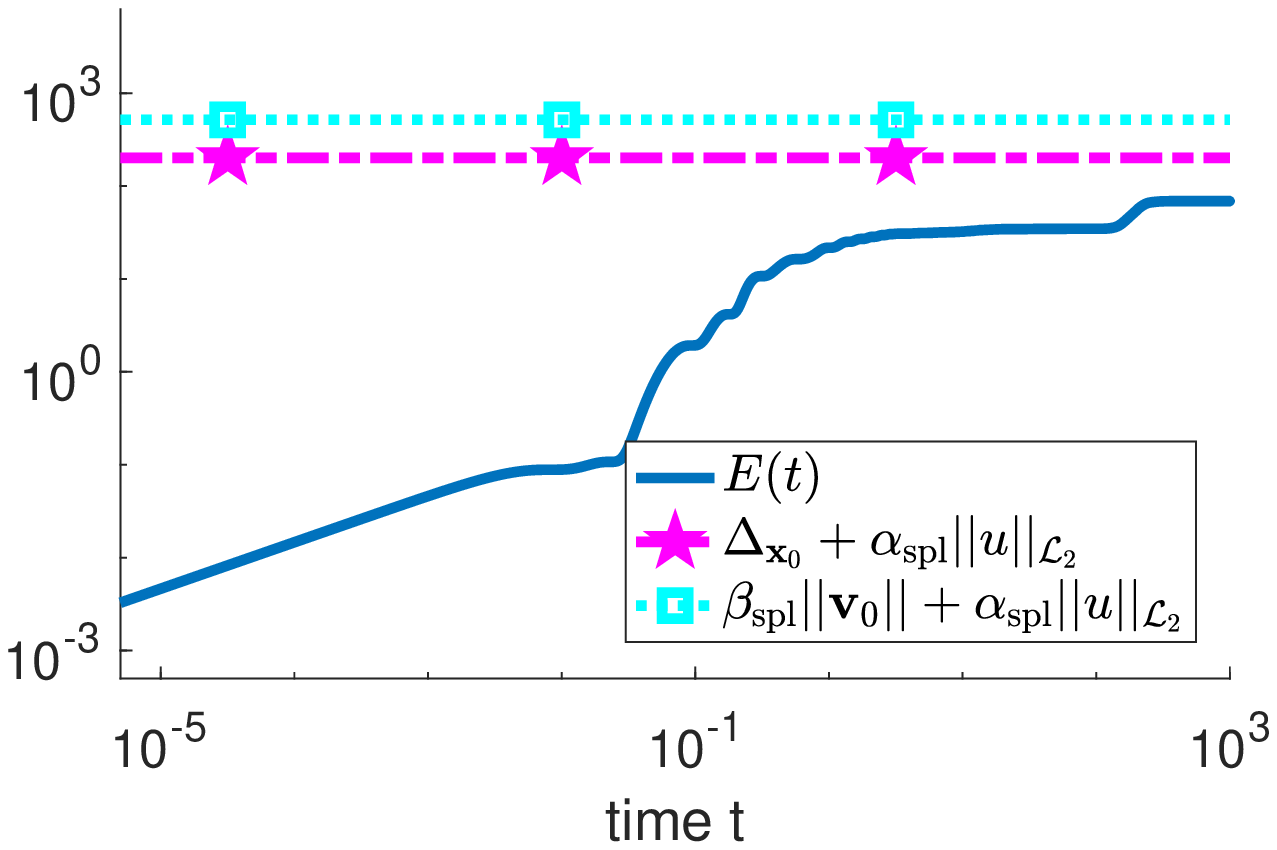}    
\caption{Beam, trained case reduced with \acro{BT-BT}. Cumulated error $E(t)$ versus estimators.
\label{fig:beamerr-tr-btbt}}  
                          
\includegraphics[height=4.1cm]{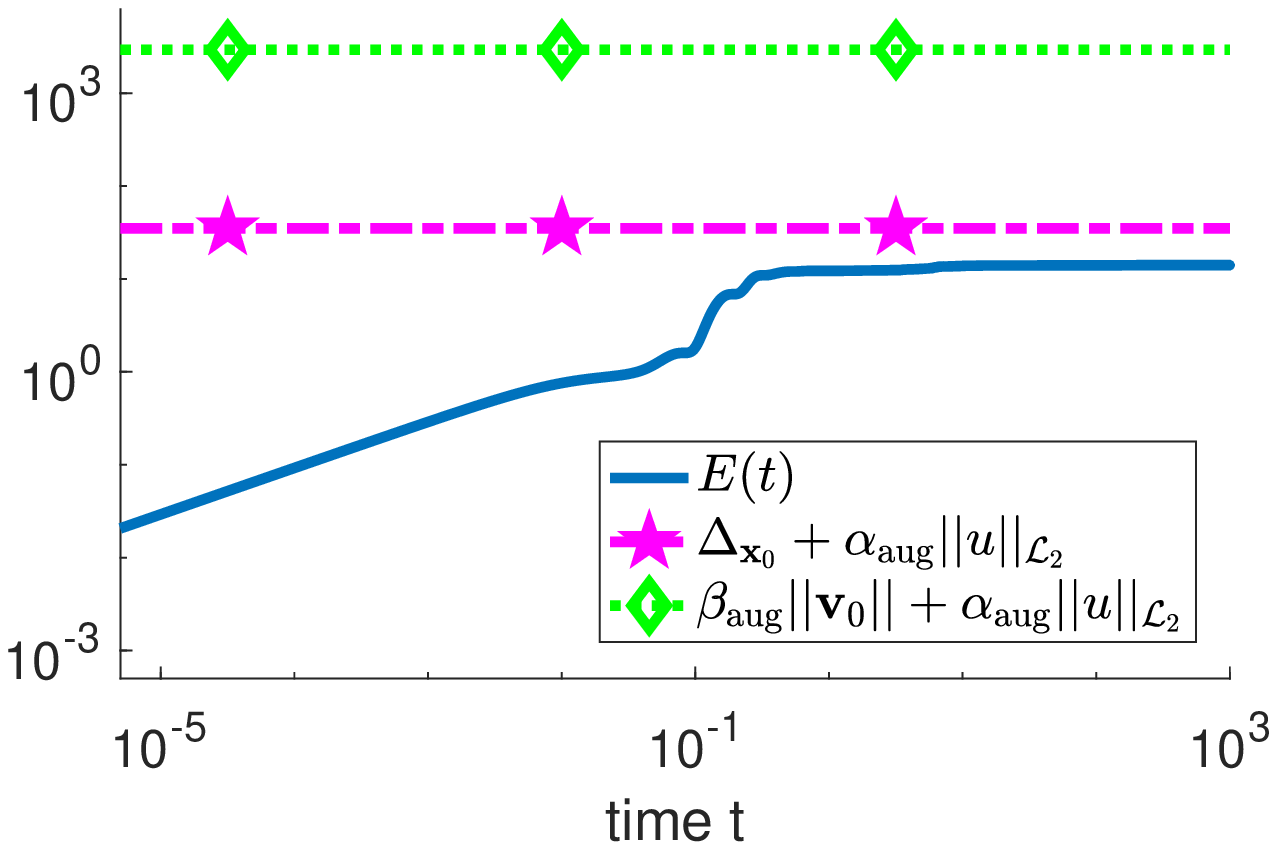}    
\caption{Beam, trained case reduced by \acro{BT-aug}. Cumulated error $E(t)$ versus estimators.
\label{fig:beamerr-tr-btaug}}                              
\end{center}                            
\end{figure}

We first consider the trained case, which is very similar to the scenario considered in \cite{art:bt-schroeder20}. The only difference is the interval in which the input is non-zero. The \acro{FOM} output is illustrated in Fig.~\ref{fig:beam-trainedsol}.
Two reduced models of dimension $n=30$ are constructed. The first one is obtained by \acro{BT-BT} with 15 degrees of freedom for the controlled and uncontrolled part each, the other by \acro{BT-aug}. Comparing Fig.~\ref{fig:beamerr-tr-btbt} and Fig.~\ref{fig:beamerr-tr-btaug}, it is seen that \acro{BT-aug} yields errors of an order lower in this parameter setting, but the error estimator for \acro{BT-BT} \eqref{eq:bound-BTBT} shows significantly less overestimation than the one for \acro{BT-aug} \eqref{eq:bound-BTaug}. Without any doubt, our proposed estimator $\Delta_{\bv x_0}$  is the most effective in predicting the error. When only considering the uncontrolled part, the difference in the effectivity becomes even more evident. Table~\ref{tab:beamerr-tr} shows that our estimator predicts the three leading digits correctly, whereas the estimate \acro{BT-BT} from literature overestimates the error by one order and the one for \acro{BT-aug} even by two orders.

Next, the not trained case is considered. As restriction \eqref{eq:x0restriction}  is violated for this case, the bounds \eqref{eq:bound-BTaug}-\eqref{eq:bound-BTBT}  are no longer applicable. Further, as this is an uncontrolled problem,  \acro{BT-BT} and \acro{BT-aug} yield the same \acro{ROM}, which can also be interpreted  as reducing the uncontrolled system by \acro{BT}. We compare this \acro{ROM} with the ones obtained using \acro{IRKA} and \acro{ISRK} instead. The errors $E(T)$ for varying order $n$ are stated in Table~\ref{tab:beam-nt} (upper entries). The overall fidelity of the \acro{IRKA} models is worst, only for the largest dimension $n=30$ the quality of the \acro{IRKA} and the \acro{BT} model are comparable. \acro{ISRK} performs similar to \acro{BT}, the highest order model with $n=30$ is best for \acro{ISRK}, which we would consider as the most favorable choice here. Our error estimator predicts all errors with a very high precision, the relative distance between estimation and error is  about a thousandth at worst, see Table~\ref{tab:beam-nt} (lower entries).

\begin{table}[tb]
{\small
\renewcommand{\arraystretch}{2.0}
 \begin{tabular}{  l |  c |  c }
        & \acro{BT-BT}  & \acro{BT-aug} \cr \hline \hline
       Error $E_{x_0}(T)$ & $3.47 \cdot 10^{1}$ & $1.40 \cdot 10^{1}$ \cr  \hline
   Proposed bound $\Delta_{\bv x_0}$ &  $ 3.47 \cdot 10^{1}$ & $1.40 \cdot 10^{1}$ \cr  \hline
     Bound literature 
     &  \makecell{$3.51 \cdot 10^{2}$\\($=\beta_{\mathrm{spl}} ||\bv v_0||)$} & \makecell{$2.91 \cdot 10^{3}$\\($=\beta_{\mathrm{aug}} ||\bv v_0|| $)}
       \end{tabular}
\renewcommand{\arraystretch}{1}       
}
\caption{Beam, trained case. Errors and estimators for uncontrolled part.
\label{tab:beamerr-tr}}
\end{table}     

{\small
\begin{table}[tb]
  \begin{tabular}{  c |  l |  l |  l |  l }
   $n$ & 12 & 18 & 24 & 30 \cr \hline \hline
     \multirow{2}{*}{\acro{BT}}   &\sepsymb{1.4 $\cdot$ 10$^{-1}$} &\sepsymb{1.8 $\cdot$ 10$^{\text{-}2}$ } &\sepsymb{3.7 $\cdot$ 10$^{\text{-}2}$ } &\sepsymb{1.1 $\cdot$ 10$^{\text{-}2}$ } \\[-0.3em]
     &5.3 $\cdot$ 10$^{\text{-}6}$ &1.3 $\cdot$ 10$^{\text{-}5}$ &2.4 $\cdot$ 10$^{\text{-}5}$ &5.5 $\cdot$ 10$^{\text{-}5}$ \cr \hline
    \multirow{2}{*}{\acro{IRKA}}   &\sepsymb{2.7\hspace{2.9em}}&\sepsymb{2.7\hspace{2.9em}}&\sepsymb{2.7\hspace{2.9em}}&\sepsymb{1.3 $\cdot$ 10$^{\text{-}2}$ } \\[-0.3em]
     &1.6 $\cdot$ 10$^{\text{-}6}$ &1.6 $\cdot$ 10$^{\text{-}6}$ &1.6 $\cdot$ 10$^{\text{-}6}$ &1.2 $\cdot$ 10$^{\text{-}4}$ \cr \hline
     \multirow{2}{*}{\acro{ISRK}}   &\sepsymb{2.7\hspace{2.9em}}&\sepsymb{2.8 $\cdot$ 10$^{\text{-}2}$ } &\sepsymb{2.8 $\cdot$ 10$^{\text{-}2}$ } &\sepsymb{6.5 $\cdot$ 10$^{\text{-}3}$ } \\[-0.3em]
    &1.5 $\cdot$ 10$^{\text{-}6}$ &1.3 $\cdot$ 10$^{\text{-}5}$ &7.8 $\cdot$ 10$^{\text{-}5}$ &6.2 $\cdot$ 10$^{\text{-}5}$
  \end{tabular}
   \caption{Beam, not trained. Error $E(T)$ and its relative difference to estimator $|\Delta_{\bv{x}_0}-E(T)| /E(T)$ {(separated by dashed line) for varying order $n$}.
\label{tab:beam-nt}}    
\end{table}   
}

\subsection{Convection-diffusion equation} \label{subsec:num-convdiff}
	\newcommand{\pdeZ}{z}
	\newcommand{\pdeXi}{\xi}	
	We consider a convection-diffusion equation on the unit square $\Omega= [0,1] \times [0,1]$. For $0<t\leq T$, $T=1$, and $\pdeXi= (\pdeXi_1,\pdeXi_2) \in \Omega$, the dynamics is described by
\begin{align*}
	\partial_t \pdeZ(t,\pdeXi) &= (\partial_{\pdeXi_1 \xi_1}+\partial_{\pdeXi_2 \pdeXi_2}) \pdeZ(t,\pdeXi) + \frac{1}{2} \pdeXi_1^2 \pdeXi_2 \partial_{\pdeXi_1} \pdeZ(t,\pdeXi).
	\end{align*}
The initial and boundary conditions are chosen as
\begin{align*}
	\pdeZ(0,\pdeXi) = x_0(\pdeXi), \, \pdeXi \in \Omega, \quad \text{and} \quad  \pdeZ(t,\pdeXi)_{\pdeXi\in |\partial \Omega} &= 0, \, t> 0.
\end{align*}
The output $\bv y(t) =[y_1(t),\ldots, y_9(t)]^T$ is defined by
\begin{align*}
	y_\ell(t) =& \int_{K_\ell} \pdeZ(t,\pdeXi) d\pdeXi, \qquad \quad \ell =1, \ldots,9, \\
	& K_\ell = \left[ \frac{9}{20} \, , \, \frac{11}{20}\right] \times \left[ \frac{\ell}{10} - \frac{1}{50} \, , \,  \frac{\ell}{10} + \frac{1}{50} \right].
\end{align*}
For the training of the \acro{ROM}, the initial values 
\begin{align*}
\tilde{x}_0^{\mu}(\pdeXi) = & \, \pdeXi_1^{\frac{1}{4}} \pdeXi_2^{\frac{1}{\mu}} (1- \pdeXi_1) (1- \pdeXi_2) \\
& \left[ \cos{\left( 10 \left( \pdeXi_2+ \frac{\mu}{5} \right)^3 \right)}  + e^{\pdeXi_1^2 \frac{\mu}{1+\pdeXi_1 \pdeXi_2}} \right]
\end{align*}
with $\mu = 2+ k/20$ for $k = 0,\ldots,20$ are used. The training matrix after space discretization has the dimension $\bt X_0 \in \R^{N,21}$. The simulation setup uses
 $x_0= \tilde{x}_0^{\mu}$ with $\mu = 3$, which lies in the column-span of $\bt X_0$.
The space discretization is done with standard finite differences on the uniform mesh with $\tilde{n}= 150$ inner points in each direction, which yields a \acro{FOM} of dimension $N= 22\,500$ with no input ($\bt B= \bt 0$) and $\bt C \in \R^{9,N}$.
For this uncontrolled large-scale problem, we examine reduced models obtained by \acro{BT}, \acro{IRKA} and \acro{ISRK}, respectively. The Gramians are determined using the sparse solvers from the \texttt{MESS} toolbox. The reduction results shown in Table~\ref{tab:convdiff} suggest that all three methods yield comparable and overall satisfactory results in terms of fidelity. Again, we find \acro{ISRK} as the most favorable choice, because the fixed point iteration underlying the construction of the \acro{ROM}s converged about twice as fast as the one of \acro{IRKA}. The latter is in correspondence with the theoretical stability guarantees that \acro{ISRK} yields. As seen in Table~\ref{tab:convdiff} (lower entries), our estimator is very effective in predicting the errors, apart from the presence of low-rank approximation errors in the Gramians. The largest relative difference to the error is observed for \acro{ISRK} with $n=30$, and there it is still only under five-tenth of a percent.

\begin{figure}[tb]
\begin{center}
\includegraphics[height=4.0cm, width=6.0cm]{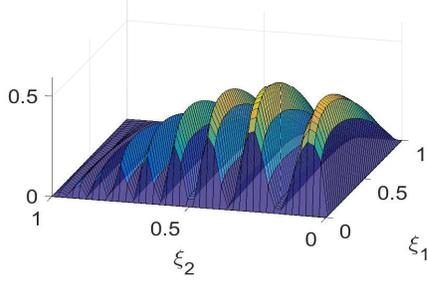}    
\caption{Convection-diffusion equation. Visualization of initial state $x_0(\pdeXi)$ used in the simulation.}  
\label{fig:convdiff-x0}                                                           
\end{center}                            
\end{figure}

{\small
\begin{table}[tb]
  \begin{tabular}{  c |  l |  l |  l |  l }
   $n$ & 12 & 18 & 24 & 30 \cr \hline \hline
     \multirow{2}{*}{\acro{BT}}  &\sepsymb{3.3 $\cdot$ 10$^{\text{-}6}$ } & \sepsymb{2.9 $\cdot$ 10$^{\text{-}7}$} & \sepsymb{5.3 $\cdot$ 10$^{\text{-}8}$ } & \sepsymb{4.2 $\cdot$ 10$^{\text{-}9}$ } \\[-0.3em]
    &1.6 $\cdot$ 10$^{\text{-}6}$ &2.2 $\cdot$ 10$^{\text{-}6}$ &8.9 $\cdot$ 10$^{\text{-}6}$ &4.5 $\cdot$ 10$^{\text{-}4}$ \cr \hline
     \multirow{2}{*}{\acro{IRKA}}  &\sepsymb{2.0 $\cdot$ 10$^{\text{-}6}$ } &\sepsymb{2.7 $\cdot$ 10$^{\text{-}7}$ } &\sepsymb{3.9 $\cdot$ 10$^{\text{-}8}$ } &\sepsymb{5.3 $\cdot$ 10$^{\text{-}9}$ } \\[-0.3em]
    &1.6 $\cdot$ 10$^{\text{-}6}$ &1.6 $\cdot$ 10$^{\text{-}6}$ &8.2 $\cdot$ 10$^{\text{-}6}$ &9.0 $\cdot$ 10$^{\text{-}4}$ \cr \hline
     \multirow{2}{*}{\acro{ISRK}}  &\sepsymb{2.3 $\cdot$ 10$^{\text{-}6}$ } &\sepsymb{2.8 $\cdot$ 10$^{\text{-}7}$ } &\sepsymb{5.4 $\cdot$ 10$^{\text{-}8}$ } &\sepsymb{4.3 $\cdot$ 10$^{\text{-}9}$ } \\[-0.3em]
    &1.6 $\cdot$ 10$^{\text{-}6}$ &2.5 $\cdot$ 10$^{\text{-}6}$ &7.5 $\cdot$ 10$^{\text{-}6}$ &4.6 $\cdot$ 10$^{\text{-}3}$
  \end{tabular}
 \caption{Convection-diffusion equation. Error $E(T)$ and its relative difference to estimator $|\Delta_{\bv{x}_0}-E(T)| /E(T)$ {(separated by dashed line) for varying order $n$}.
\label{tab:convdiff}}    
\end{table}   
}
\section*{Conclusion}

We showed that the reduction error related to inhomogeneous initial conditions can be estimated very effectively using the observability Gramian of the error system. Our results are applicable for the reduction of linear time invariant systems that yield asymptotically stable models and can be efficiently implemented given the observability Gramian of the full order model is available. As such, they could be a building block for certified model reduction with inhomogeneous initial conditions. Possible related topics for future research are e.g., the application in the reduction of switched systems or the adaption of our approach to time-limited model reduction.

\begin{ack}                               
The author acknowledges the support of the DFG research training group 2126 on algorithmic optimization.
\end{ack}

\bibliographystyle{plain}        
\bibliography{BtIc}           

\begin{thebibliography}{10}

\bibitem{book:antoulas2005}
A.~Antoulas.
\newblock {\em Approximation of Large-Scale Dynamical Systems}, volume~6 of
  {\em Adv. Des. Control}.
\newblock {SIAM}, Philadelphia, PA, 2005.

\bibitem{morAntBG10}
A.~C. Antoulas, C.~A. Beattie, and S.~Gugercin.
\newblock Interpolatory model reduction of large-scale dynamical systems.
\newblock In Javad Mohammadpour and Karolos~M. Grigoriadis, editors, {\em
  Efficient Modeling and Control of Large-Scale Systems}, pages 3--58. Springer
  US, 2010.

\bibitem{art:morAntSZ02}
A.~C. Antoulas, D.~C. Sorensen, and Y.~Zhou.
\newblock On the decay rate of {H}ankel singular values and related issues.
\newblock {\em Systems Control Lett.}, 46(5):323--342, 2002.

\bibitem{morBeaGM17}
C.~Beattie, S.~Gugercin, and V.~Mehrmann.
\newblock Model reduction for systems with inhomogeneous initial conditions.
\newblock {\em Systems Control Lett.}, 99:99--106, 2017.

\bibitem{book:dimred2003}
P.~Benner, V.~Mehrmann, and D.~C. Sorensen, editors.
\newblock {\em Dimension Reduction of Large-Scale Systems}.
\newblock Lecture Notes in Computational Science and Engineering. Springer, 1
  edition, 2005.

\bibitem{art:Daraghmeh2019}
A.~Daraghmeh, C.~Hartmann, and N.~Qatanani.
\newblock Balanced model reduction of linear systems with nonzero initial
  conditions: singular perturbation approximation.
\newblock {\em Appl. Math. Comput.}, 353:295--307, 2019.

\bibitem{book:LinearRobustControl}
M.~Green and D.~J.~N. Limebeer.
\newblock {\em Linear Robust Control}.
\newblock Prentice-Hall, Inc., New York, 1994.

\bibitem{art:isrk-Gugercin08}
S.~Gugercin.
\newblock An iterative {SVD}-krylov based method for model reduction of
  large-scale dynamical systems.
\newblock {\em Linear Algebra Appl.}, 428(8--9):1964--1986, 2008.

\bibitem{art:morHeiRA11}
M.~Heinkenschloss, T.~Reis, and A.~C. Antoulas.
\newblock Balanced truncation model reduction for systems with inhomogeneous
  initial conditions.
\newblock {\em Automatica}, 47(3):559--564, 2011.

\bibitem{code:bls22-btinhom}
B.~Liljegren-Sailer.
\newblock Code for the paper {'Effective error estimation for model reduction
  with inhomogeneous initial conditions'}.
\newblock https://doi.org/10.5281/zenodo.5863661, 2022.

\bibitem{morMoo81}
B.~C. Moore.
\newblock Principal component analysis in linear systems: controllability,
  observability, and model reduction.
\newblock {\em {IEEE} Trans. Autom. Control}, AC--26(1):17--32, 1981.

\bibitem{inproc:PoussotVassal2011}
C.~Poussot-Vassal.
\newblock An iterative {SVD}-tangential interpolation method for medium-scale
  {MIMO} systems approximation with application on flexible aircraft.
\newblock In {\em {IEEE} Conference on Decision and Control and European
  Control Conference}. {IEEE}, 2011.

\bibitem{SaaKB-mmess}
J.~Saak, M.~K\"{o}hler, and P.~Benner.
\newblock {M-M.E.S.S.-2.1} -- {T}he {M}atrix {E}quations {S}parse {S}olvers
  library.
\newblock 10.5281/zenodo.4719688, 2021.

\bibitem{art:bt-schroeder20}
C.~Schr{\"o}der and M~Voigt.
\newblock Balanced truncation model reduction with a priori error bounds for
  {LTI} systems with nonzero initial value.
\newblock arXiv e-prints 2006.02495, 2020.

\bibitem{art:Simoncini2016}
V.~Simoncini.
\newblock Computational methods for linear matrix equations.
\newblock {\em {SIAM} Rev.}, 58(3):377--441, 2016.

\bibitem{inproc:Vasilyev05}
D.~Vasilyev and J.~White.
\newblock A more reliable reduction algorithm for behavioral model extraction.
\newblock In {\em {ICCAD}-2005. {IEEE}/{ACM} International Conference on
  Computer-Aided Design}. {IEEE}, 2005.

\bibitem{code:voigt20-btinhom}
M.~Voigt.
\newblock {BT\_INHOM -- MATLAB routines for balanced truncation model reduction
  for systems with nonzero initial value}.
\newblock https://doi.org/10.5281/zenodo.3875468, 2020.

\bibitem{art:TSIA-Xu11}
Y.~Xu and T.~Zeng.
\newblock Optimal {H2} model reduction for large scale {MIMO} systems via
  tangential interpolation.
\newblock {\em Int. J. Numer. Anal. Model.}, 8(1):174--188, 2011.

\end{thebibliography}



\end{document}